\documentclass[11pt]{article}
\usepackage[T1]{fontenc}
\usepackage{lmodern}
\usepackage{amsmath,amssymb,amsthm,mathtools}
\usepackage{microtype}
\usepackage{xurl}
\usepackage{xcolor}
\usepackage{enumitem}
\usepackage{hyperref}
\hypersetup{
  colorlinks=true,
  linkcolor=blue!55!black,
  citecolor=blue!55!black,
  urlcolor=blue!55!black,
  pdftitle={Partition Principle without Choice via Symmetric Iterations and Sheaf--Toposes: Retraction of the Proposed Construction},
  pdfauthor={Frank Gilson}
}

\newcommand{\ZF}{\mathsf{ZF}}
\newcommand{\IZF}{\mathsf{IZF}}
\newcommand{\AC}{\mathsf{AC}}
\newcommand{\DC}{\mathsf{DC}}
\newcommand{\PP}{\mathsf{PP}}
\newcommand{\Sh}{\mathsf{Sh}}
\newcommand{\NN}{\mathbb N}
\newcommand{\onto}{\twoheadrightarrow}
\newcommand{\into}{\hookrightarrow}

\theoremstyle{plain}
\newtheorem{proposition}{Proposition}[section]

\theoremstyle{definition}
\newtheorem{definition}[proposition]{Definition}

\theoremstyle{remark}

\title{Partition Principle without Choice via Symmetric Iterations and Sheaf--Toposes\\
\large Retraction of the Proposed Construction}
\author{Frank Gilson}
\date{\today}

\begin{document}
\maketitle

\begin{abstract}
An earlier form of this manuscript proposed that a transformation-groupoid
topos associated with a free action of a nontrivial finite group on Cantor
space contains an internal model of
\(\IZF+\PP+\neg\AC\), with a classical upgrade after
double-negation sheafification. That construction is not valid.
The proposed local embedding argument does not establish the required
property for general epimorphisms; the finite-fiber class of small maps does
not support an \(\IZF\)-universe because it does not contain the natural
numbers object; the quotient map used to witness failure of choice is not a
nonsplitting epimorphism in the asserted form; and double-negation
sheafification does not repair these defects. The symmetric-iteration
appendix also used an invalid same-condition equivariance inference.
Accordingly, this replacement withdraws all model-existence, preservation,
and independence claims. We retain the conditional categorical observation
that a genuinely matching local family of embeddings descends to a global
embedding, and record the precise obstructions so that the withdrawn
construction is not cited as a solution of the Partition Principle problem.
\end{abstract}



\section{Revision notice and status}

The Partition Principle is
\[
  \PP:\qquad
  \forall A\,\forall B\,
  \bigl((\exists f:A\onto B)\longrightarrow
        (\exists g:B\into A)\bigr).
\]
The implication \(\AC\Rightarrow\PP\) is immediate. The converse over
\(\ZF\) remains a longstanding open problem; see
\cite{JechAC,HowardRubin}. Nothing in this paper resolves that problem.

The earlier manuscript claimed an internal set universe \(V\) in
\[
  \mathcal E=\Sh(H\ltimes 2^{\NN}),
\]
where \(H\) is a nontrivial finite group acting freely on Cantor space,
such that
\[
  V\models\IZF+\PP+\neg\AC,
\]
and claimed a classical \(\ZF\)-model internally after
double-negation sheafification. It also presented symmetric-extension
templates purporting to preserve failures of \(\PP\) and \(\DC\).

\medskip
\noindent\textbf{All of those claims are withdrawn.}
They depend on several independent failures described below. In particular,
the former main theorem, its preview, the asserted internal
\(\PP\)-proposition, the asserted internal \(\neg\AC\)-proposition, the
Booleanization upgrade, and the preservation appendix must not be used.

\section{The conditional descent observation}

The one categorical implication retained from the proposed route is
conditional and essentially formal.

\begin{definition}[Matching local embeddings]
Let \(\mathcal E\) be a Grothendieck topos. An epimorphism
\(e:A\onto B\) has the matching local embedding property if there is a
covering family \(\{U_i\to 1\}\) and monomorphisms
\[
  s_i:B\times U_i\into A\times U_i
\]
whose restrictions agree as morphisms on every overlap
\(U_i\times U_j\).
\end{definition}

\begin{proposition}
If \(e:A\onto B\) has the matching local embedding property, then there is
a global monomorphism \(s:B\into A\).
\end{proposition}

\begin{proof}
Morphisms into \(A\) satisfy descent along an effective epimorphic family.
The matching arrows \(s_i\) therefore glue to a unique arrow \(s:B\to A\).
Monicity is local for the canonical topology, so \(s\) is monic.
\end{proof}

This proposition does not supply the matching family. With matching included
in the hypothesis, the proposed \emph{local embedding property} already
contains the descent data needed to obtain the desired global monomorphism.
Thus the entire mathematical burden is proving that property for every
relevant epimorphism. The earlier manuscript did not do so.

\section{Failure of the local embedding argument}

The general-fiber proof attempted to write an equivariant etale domain as a
filtered union
\[
  A=\mathop{\mathrm{colim}}_{F} A_F
\]
of finite unions of local slices, apply a finite-fiber argument to maps
\(e_F:A_F\to B\), and pass compatible embeddings
\(B\into A_F\) to the colimit. Three problems are fatal.

\begin{enumerate}[label=(\arabic*)]
\item A finite union \(A_F\) of slices need not meet every fiber over \(B\).
Consequently \(e_F:A_F\to B\) need not be an epimorphism, so the claimed
finite-fiber lemma is inapplicable.

\item For \(F\subseteq G\), the manuscript used a supposed projection
\[
  [\operatorname{id}_{A_F},0]:
  A_F\amalg A_{G\setminus F}\longrightarrow A_F.
\]
A general topos is not pointed: there need not be any arrow
\(A_{G\setminus F}\to A_F\). Hence this projection is not available.

\item Even if stagewise embeddings existed, arbitrary choices of branches do
not form a compatible system. Exactness of filtered colimits preserves a
compatible diagram; it cannot manufacture compatibility among maps that were
chosen independently.
\end{enumerate}

The finite-fiber precursor also does not establish the stated equivariant
claim. A continuous local branch of an etale map is not automatically
equivariant. Moreover, a clopen bisection in the arrow space of an action
groupoid is not, merely by being a bisection, one of the asserted disjoint
covering objects of the terminal object. These distinctions cannot be
removed by refining the topology.

It follows that no local embedding property for all relevant epimorphisms was
proved. Therefore the descent proposition above gives no instance of
\(\PP\).

\section{Failure of the algebraic-set-theory construction}

The earlier manuscript declared a morphism small when its corresponding
etale map has finite fibers. In particular, an object \(X\) is small only
when \(X\to 1\) has finite fibers. The natural numbers object
\(\NN\to 1\) is therefore not small.

Smallness of the natural numbers object is the standard axiom
\((\mathrm{NS})\) used to interpret Infinity in an initial ZF-algebra.
To obtain \(\IZF\), the usual algebraic-set-theory theorem additionally
requires the appropriate Power Set and Separation axioms for the class of
small maps; see \cite{JoyalMoerdijk,vdBMUnified,vdBMPredicative}. Closure
under pullback, composition, finite coproducts, and a proposed collection
argument is not a substitute for these hypotheses.

Thus the finite-fiber construction can at most describe a
hereditarily-finite-style universe. It cannot be the asserted model of
\(\IZF\). References to an initial algebra do not make Infinity, Power Set,
or full Separation hold when their small-map hypotheses fail.

There was a second, independent transfer gap. The manuscript treated
arbitrary morphisms between ambient \emph{displays} as if they were fully and
conservatively identical with maps between sets of the initial algebra.
No theorem establishing that strong correspondence for the displays used in
the argument was provided. Accordingly, neither an ambient epimorphism nor
an ambient monomorphism could simply be re-labelled as the required internal
set map.

\section{Failure of the proposed witness to
\texorpdfstring{\(\neg\AC\)}{not-AC}}

Let a finite group \(H\) act freely on a compact zero-dimensional Hausdorff
space \(X\), and let
\[
  q:X\longrightarrow X/H
\]
be the quotient map. The previous argument regarded \(q\) simultaneously
as an equivariant etale morphism over \(X\) and as a covering over \(X/H\).
These are different categorical presentations. In the category of
\(H\)-equivariant etale spaces over \(X\), the displayed codomain \(X/H\)
does not naturally carry the required structure map to \(X\).

Under the quotient-space interpretation the proposed nonsplitting conclusion
is false.

\begin{proposition}
The finite covering \(q:X\to X/H\) has a continuous section.
\end{proposition}

\begin{proof}
The quotient \(X/H\) is compact, Hausdorff, and zero-dimensional. Choose a
finite open cover over which \(q\) is trivial. A finite open cover of a
compact zero-dimensional space has a finite pairwise-disjoint clopen
refinement. On every member of that refinement choose one sheet of the
trivialized covering. The corresponding local sections glue, since their
domains are disjoint and clopen, to a continuous global section of \(q\).
\end{proof}

The earlier fixed-point calculation proves only that there is no
\(H\)-equivariant map from \(X/H\), equipped with the trivial action, into
the free \(H\)-space \(X\). It does not prove that the covering object in the
quotient sheaf topos has no section. Hence \(q\) cannot witness failure of
internal choice in the proposed construction.

\section{Why Booleanization does not repair the construction}

Double-negation sheafification is left exact, so it preserves finite limits
and monomorphisms. It does not follow that it preserves the assertion that a
particular epimorphism has no section: a reflector may create new morphisms,
including sections. Nor does left exactness establish the missing local
embedding theorem.

More fundamentally, Booleanization does not turn a finite-fiber small-map
class into one containing the natural numbers object, and it does not
automatically transport an initial ZF-algebra together with all the
small-map axioms needed for \(\IZF\). A classical internal model would
require a new, complete verification in the Boolean subtopos. No such
verification is available here.

The claimed upgrade to \(\ZF+\PP+\neg\AC\) is therefore withdrawn
independently of the earlier intuitionistic claim.

\section{Failure of the symmetric-extension appendix}

The single-stage obstruction used the following invalid inference. If
\(\pi\) is a forcing automorphism, the automorphism lemma gives
\[
  p\Vdash\varphi(\dot x)
  \quad\Longleftrightarrow\quad
  \pi p\Vdash\varphi(\pi\dot x).
\]
Even when \(\pi\dot g=\dot g\), information forced by \(p\) is transported to
\(\pi p\), not back to \(p\). To obtain a contradiction one still needs, for
example, that \(\pi p=p\), or that \(p\) and \(\pi p\) have a common extension
that decides incompatible values. Transitivity of a fiber action and support
of the whole name do not supply that same-condition fact.

The iteration argument has a further mismatch. It chooses a later coordinate
outside the support and claims that its automorphisms move a value of a
function \(g:B\to A\) where \(A\) and \(B\) were introduced at an earlier
stage. Later-coordinate automorphisms normally act trivially on those old
objects, so they do not move \(g(b)\in A\).

Finally, the asserted preservation of \(\DC\) was not proved by choosing a
successor locally at each finite stage. Producing a countable branch requires
one coherent hereditarily symmetric name and control of the union of its
supports. A level-by-level list of local choices does not provide that
object.

All negative-\(\PP\), iteration-preservation, and \(\DC\)-preservation claims
from that appendix are withdrawn.

\section{Conclusion}

The proposed sheaf-topos and symmetric-iteration construction does not
produce a model of
\[
  \IZF+\PP+\neg\AC
  \qquad\text{or}\qquad
  \ZF+\PP+\neg\AC.
\]
The defects occur before any question of externalizing an internal model:
the local embedding theorem is unproved, the finite-fiber small-map universe
fails Infinity, and the proposed quotient epimorphism does not witness
failure of choice. The symmetric appendix supplies no independent rescue.

This correction is intentionally definitive. The original construction and
all consequences depending on it are superseded by the present account. The
general question whether \(\ZF+\PP+\neg\AC\) is consistent remains open.

\end{document}